\documentclass[12pt,a4paper,twoside]{article}

\newcommand{\pageformat}[6]{\setlength{\hoffset}{-1in}
                  \setlength{\voffset}{-1in}
                  \addtolength{\hoffset}{#5}
                            \addtolength{\voffset}{#6}
                            \setlength{\oddsidemargin}{#1}
                            \setlength{\evensidemargin}{#2}
                            \setlength{\textwidth}{\paperwidth}
                  \addtolength{\textwidth}{-\oddsidemargin}
                  \addtolength{\textwidth}{-\evensidemargin}
                  \addtolength{\textwidth}{-\marginparsep}
                  \addtolength{\textwidth}{-\marginparwidth}
                            \setlength{\topmargin}{#3}
                            \setlength{\textheight}{\paperheight}
                  \addtolength{\textheight}{-\topmargin}
                  \addtolength{\textheight}{-\headheight}
                  \addtolength{\textheight}{-\headsep}
                  \addtolength{\textheight}{-\footskip}
                  \addtolength{\textheight}{-#4}}
\pageformat{2cm}{3cm}{24mm}{24mm}{1pt}{0pt}

\usepackage{ifthen}
\newboolean{article}
    \setboolean{article}{true}
\newboolean{report}
\newboolean{book}
\newboolean{letter}
\newboolean{german}
\newboolean{italian}
\newboolean{nobaselinestretch}
\newboolean{nosectionappendix}
\newboolean{oldtoc}
\newboolean{nosectionequn}
\newboolean{notheorem}

\ifthenelse{\boolean{german}}{
    \usepackage{german}}{}

\usepackage[latin1]{inputenc}
\usepackage[hypertex]{hyperref}

\usepackage{amsmath}
\usepackage{amssymb}
\usepackage[mathscr]{eucal}

\ifthenelse{\boolean{notheorem}}{}{
    \usepackage{theorem}}

\ifthenelse{\boolean{nobaselinestretch}}{}{
    \renewcommand{\baselinestretch}{1.25}}

\newenvironment{env}[2]{\begin{#1}#2\end{#1}}{}
    \newcommand{\beq}[1]{\begin{env}{equation}{#1}}
    \newcommand{\beqn}[1]{\begin{env}{equation*}{#1}}
    \newcommand{\bal}[1]{\begin{env}{align}{#1}}
    \newcommand{\baln}[1]{\begin{env}{align*}{#1}}
    \newcommand{\bga}[1]{\begin{env}{gather}{#1}}
    \newcommand{\bgan}[1]{\begin{env}{gather*}{#1}}
    \newcommand{\bflal}[1]{\begin{env}{flalign}{#1}}
    \newcommand{\bflaln}[1]{\begin{env}{flalign*}{#1}}
    \newcommand{\bmu}[1]{\begin{env}{multline}{#1}}
    \newcommand{\bmun}[1]{\begin{env}{multline*}{#1}}
    \newcommand{\bsp}[1]{\begin{env}{split}{#1}}

    \newcommand{\eeq}{\end{env}}
    \newcommand{\eeqn}{\end{env}}
    \newcommand{\eal}{\end{env}}
    \newcommand{\ealn}{\end{env}}
    \newcommand{\ega}{\end{env}}
    \newcommand{\egan}{\end{env}}
    \newcommand{\eflal}{\end{env}}
    \newcommand{\eflaln}{\end{env}}
    \newcommand{\emu}{\end{env}}
    \newcommand{\emun}{\end{env}}
    \newcommand{\esp}{\end{env}}

\newcommand{\lf}{\vspace{2ex}}

\renewcommand{\bf}[1]{\textbf{#1}}
\renewcommand{\it}[1]{\textit{#1}}

\renewcommand{\sf}[1]{\textsf{#1}}

\renewcommand{\tt}[1]{\texttt{#1}}
\newcommand{\hl}[1]{\bf{\it{#1}}}

\newcommand{\mbf}[1]{\mathbf{#1}}
\newcommand{\msf}[1]{\text{\small$\sf{#1}$}}

\newcommand{\cmc}[1]{\mathcal{#1}}
\newcommand{\eus}[1]{\mathscr{#1}}
\newcommand{\euf}[1]{\mathfrak{#1}}
\newcommand{\bb}[1]{\mathbb{#1}}

\newcommand{\nbd}[1]{$#1$\nobreakdash--}
\newcommand{\ol}[1]{\overline{#1}}
\newcommand{\ul}[1]{\underline{#1}}
\newcommand{\wt}[1]{\widetilde{#1}}
\newcommand{\wh}[1]{\widehat{#1}}

\newcommand{\vt}{\vartheta}

\newcommand{\Om}{\Omega}

\newcommand{\bfam}[1]{\bigl(#1\bigr)}

\newcommand{\AB}[1]{\langle#1\rangle}

\newcommand{\BAB}[1]{\Bigl\langle#1\Bigr\rangle}
\newcommand{\CB}[1]{\{#1\}}
\newcommand{\bCB}[1]{\bigl\{#1\bigr\}}
\newcommand{\BCB}[1]{\Bigl\{#1\Bigr\}}

\newcommand{\LO}[1]{(#1]}

\newcommand{\set}[2][]{
    \ifthenelse{\equal{#1}{}}{
        \CB{#2}}{
        \CB{#1~|~#2}}}
\newcommand{\bset}[2][]{
    \ifthenelse{\equal{#1}{}}{
        \bCB{#2}}{
        \bCB{#1~|~#2}}}
\newcommand{\Bset}[2][]{
    \ifthenelse{\equal{#1}{}}{
        \BCB{#2}}{
        \BCB{#1~\big|~#2}}}
\newcommand{\zero}{\CB{0}}

\DeclareMathOperator{\ls}{\normalfont\msf{span}}
\DeclareMathOperator{\cls}{\ol{\ls}}
\DeclareMathOperator*{\limind}{lim\,ind}

\DeclareMathOperator{\id}{\normalfont\msf{id}}

\renewcommand{\ker}{\operatorname{\msf{ker}}}

\newcommand{\C}{\bb{C}}

\newcommand{\bJ}{\bb{J}}

\newcommand{\N}{\bb{N}}

\newcommand{\R}{\bb{R}}
\newcommand{\bS}{\bb{S}}

\newcommand{\cA}{\cmc{A}}
\newcommand{\cB}{\cmc{B}}

\newcommand{\sB}{\eus{B}}

\newcommand{\sE}{\eus{E}}

\newcommand{\sK}{\eus{K}}

\newcommand{\es}{\mathbbm{s}}
\newcommand{\et}{\mathbbm{t}}

\newcommand{\U}{\mbf{1}}

\ifthenelse{\boolean{nosectionequn}}{}{
    \numberwithin{equation}{section}
    }

\ifthenelse{\boolean{article}\or\boolean{letter}\or\boolean{nosectionequn}}{
    \setboolean{nosectionappendix}{true}}{}
\ifthenelse{\boolean{nosectionappendix}}{}{
    \renewcommand{\appendix}{
        \chapter*{\appendixname}
        \addcontentsline{toc}{chapter}{\appendixname}
        \renewcommand{\thesection}{\Alph{section}}
        \setcounter{section}{0}}}
   
\ifthenelse{\boolean{report}\or\boolean{book}}{
    }{}

\ifthenelse{\boolean{notheorem}}{}{
        \newcommand{\mnname}{Mathematical note.}
        \newcommand{\enname}{End of the note.}
        \newcommand{\definame}{Definition.}
        \newcommand{\propname}{Proposition.}
        \newcommand{\lemname}{Lemma.}
        \newcommand{\exname}{Example.}
        \newcommand{\exername}{Exercise.}
        \newcommand{\remname}{Remark.}
        \newcommand{\obname}{Observation.}
        \newcommand{\thmname}{Theorem.}
        \newcommand{\corname}{Corollary.}
        \newcommand{\proofname}{Proof.}
    \ifthenelse{\boolean{german}}{
        \renewcommand{\mnname}{Mathematische Notiz.}
        \renewcommand{\enname}{Ende der Notiz.}
        \renewcommand{\exname}{Beispiel.}
        \renewcommand{\exername}{Übung.}
        \renewcommand{\remname}{Bemerkung.}
        \renewcommand{\obname}{Beobachtung.}
        \renewcommand{\thmname}{Satz.}
        \renewcommand{\corname}{Korollar.}
        \renewcommand{\proofname}{Beweis.}}{}
    \ifthenelse{\boolean{italian}}{
        \renewcommand{\mnname}{Nota matematica.}
        \renewcommand{\enname}{Fina della nota.}
        \renewcommand{\definame}{Definizione.}
        \renewcommand{\propname}{Proposizione.}
        \renewcommand{\exname}{Esempio.}
        \renewcommand{\exername}{Esercizio.}
        \renewcommand{\remname}{Nota.}
        \renewcommand{\obname}{Osservazione.}
        \renewcommand{\thmname}{Teorema.}
        \renewcommand{\corname}{Corollario.}
        \renewcommand{\proofname}{Dimostrazione.}

       \renewcommand{\appendixname}{Appendice}

       }{}
    \theoremheaderfont{\normalfont\bfseries}
    \theoremstyle{change}
        \theorembodyfont{\rmfamily}
            \newtheorem{emp}{}[section]
                \newcommand{\bemp}[1][]{
                    \begin{emp}\hskip-\labelsep\bf{#1}\hskip\labelsep}
                \newcommand{\eemp}{\end{emp}}
\newtheorem{itemp}[emp]{}
                \newcommand{\bitemp}[1][]{
                    \begin{itemp}\hskip-\labelsep\bf{#1}\hskip\labelsep\normalfont\itshape}
                \newcommand{\eitemp}{\end{itemp}}
            \newtheorem{mn}[emp]{\mnname}
                \newcommand{\bnm}{\begin{mn}~\begin{quotation}\renewcommand{\baselinestretch}{1}\small\noindent\ignorespaces}
                \newcommand{\enm}{\end{quotation}\hfill\bf{\enname}\end{mn}}
            \newtheorem{ex}[emp]{\exname}
                \newcommand{\bex}{\begin{ex}}
                \newcommand{\eex}{\end{ex}}
            \newtheorem{exer}[emp]{\exername}
                \newcommand{\bexer}{\begin{exer}}
                \newcommand{\eexer}{\end{exer}}
            \newtheorem{defi}[emp]{\definame}
                \newcommand{\bdefi}{\begin{defi}}
                \newcommand{\edefi}{\end{defi}}
            \newtheorem{rem}[emp]{\remname}
                \newcommand{\brem}{\begin{rem}}
                \newcommand{\erem}{\end{rem}}
            \newtheorem{ob}[emp]{\obname}
                \newcommand{\bob}{\begin{ob}}
                \newcommand{\eob}{\end{ob}}
        \theorembodyfont{\normalfont\itshape}
            \newtheorem{thm}[emp]{\thmname}
                \newcommand{\bthm}{\begin{thm}}
                \newcommand{\ethm}{\end{thm}}
            \newtheorem{prop}[emp]{\propname}
                \newcommand{\bprop}{\begin{prop}}
                \newcommand{\eprop}{\end{prop}}
            \newtheorem{cor}[emp]{\corname}
                \newcommand{\bcor}{\begin{cor}}
                \newcommand{\ecor}{\end{cor}}
            \newtheorem{lem}[emp]{\lemname}
                \newcommand{\blem}{\begin{lem}}
                \newcommand{\elem}{\end{lem}}
\newenvironment{empn}[1]{\lf\noindent\bf{#1}\ignorespaces\hskip\labelsep}{\lf}
		\newcommand{\bempn}[1]{\begin{empn}{#1}}
		\newcommand{\eempn}{\end{empn}}
		\newcommand{\bitempn}[1]{\begin{empn}{#1}\normalfont\itshape}
		\newcommand{\eitempn}{\end{empn}}
                \newcommand{\bnmn}{\begin{empn}{\mnname}~\begin{quotation}\renewcommand{\baselinestretch}{1}\small\noindent\ignorespaces}
                \newcommand{\enmn}{\end{quotation}\hfill\bf{\enname}\end{empn}}
		\newcommand{\bexn}{\begin{empn}{\exname}}
		\newcommand{\eexn}{\end{empn}}
		\newcommand{\bexern}{\begin{empn}{\exername}}
		\newcommand{\eexern}{\end{empn}}
		\newcommand{\bdefin}{\begin{empn}{\definame}}
		\newcommand{\edefin}{\end{empn}}
		\newcommand{\bremn}{\begin{empn}{\remname}}
		\newcommand{\eremn}{\end{empn}}
		\newcommand{\bobn}{\begin{empn}{\obname}}
		\newcommand{\eobn}{\end{empn}}

\newcommand{\qedsymbol}{~\rule[-0.35mm]{2mm}{2mm}}
    \newcounter{proof}[emp]
    \newenvironment{Proof}[1]{
        \vspace{1ex}
        \renewcommand{\item}[1][\stepcounter{proof}(\roman{proof})]%
            {##1\hskip\labelsep}
        \noindent\textsc{#1\hskip\labelsep}}{
        \nolinebreak\qedsymbol}
    \newcommand{\proof}[1][\proofname]{
        \begin{Proof}{#1}\ignorespaces}
    \newcommand{\qed}{\end{Proof}}
    \newcommand{\noqed}{
        \renewcommand{\qedsymbol}{}
        \end{Proof}}}
    \ifthenelse{\boolean{italian}}{
        \renewcommand{\proofname}{Dimostrazione.}}{}

\usepackage[varg]{txfonts}

\renewcommand{\thefootnote}{[\alph{footnote}]}

\newcommand{\cvN}{\euf{cvN}}
\newcommand{\ecC}{\euf{cC}}

\usepackage{bbm}

\begin{document}

\title{Isometric Dilations of\\Representations of Product Systems\\via Commutants\renewcommand{\thefootnote}{}\thanks{2000 AMS-Subject classification: 46L55; 46L08; 46L53; 60J25}}
\author{}
\author{
~\\
Michael Skeide\thanks{This work is supported by research funds of University of Molise and Italian MIUR (PRIN 2005).}\\[1ex]
{\small\itshape Universit\`a\ degli Studi del Molise}\\
{\small\itshape Dipartimento S.E.G.e S.}\\
{\small\itshape Via de Sanctis}\\
{\small\itshape 86100 Campobasso, Italy}\\
{\small{\itshape E-mail: \tt{skeide@unimol.it}}}\\
{\small{\itshape Homepage: \tt{http://www.math.tu-cottbus.de/INSTITUT/lswas/\_skeide.html}}}\\
\\
}
\date{February 2006, revised November 2006}

{
\renewcommand{\baselinestretch}{1}
\maketitle




\begin{abstract}
\noindent
We construct a weak dilation of a not necessarily unital CP-semigroup to an \nbd{E}semigroup acting on the adjointable operators of a Hilbert module with a unit vector. We construct the dilation in such a way that the dilating \nbd{E}semigroup has a pre-assigned product system. Then, making use of the commutant of von Neumann correspondences, we apply the dilation theorem to proof that covariant representations of product systems admit isometric dilations.
\end{abstract}

}



\section{Introduction}\label{intro}

Let $\bS=\R_+$ or $\bS=\N_0$. Our scope is the proof of the following theorem on existence of isometric dilations of covariant representations of product systems, a problem suggested by one of the authors of \cite{MuSo02} on the occasion of a meeting in Bangalore in December 2005.

\bthm\label{CRdilthm}
Let $F^\odot=\bfam{F_t}_{t\in\bS}$ be a product system of correspondences over a \nbd{C^*}algebra (a \nbd{W^*}algebra) $M$. Let $\sigma^\odot=\bfam{\sigma_t}_{t\in\bS}$ be a (normal) covariant representation of $F^\odot$ on a Hilbert space $G$. Then there exists a (normal) isometric dilation $\tau^\odot=\bfam{\tau_t}_{t\in\bS}$ of $\sigma^\odot$.
\ethm

We recall briefly the definitions of the three notions that occur in the theorem (in a variant adapted to the theorem). By $\odot$ we indicate the \it{internal} tensor product of two correspondences.\footnote{Recall that a correspondence has always a nondegenerate left action.}  A \hl{product system} is a family $F^\odot=\bfam{F_t}_{t\in\bS}$ of correspondences $F_t$ over a \nbd{C^*}algebra $M$ with $E_0=M$, together with a family of bilinear unitaries $u_{t,s}\colon F_t\odot F_s\rightarrow F_{t+s}$ that compose associatively, and such that the identifications $u_{0,t}(F_0\odot F_t)=F_t=u_{t,0}(F_t\odot F_0)$ are the canonical ones. In the \nbd{W^*}case the tensor product is that of \nbd{W^*}correspondences. Often, we will suppress the mappings $u_{t,s}$ and simply write $F_t\odot F_s=F_{t+s}$. A \hl{covariant representation} of a product system $F^\odot$ on a Hilbert space $G$ is a family $\sigma^\odot=\bfam{\sigma_t}_{t\in\bS}$ of linear CC-maps $\sigma_t\colon F_t\rightarrow\sB(G)$ such that $\sigma_0$ is a representation of $M$ (not necessarily nondegenerate but of course $*$), all $\sigma_t$ are $M$--$M$--linear (that is, $\sigma_t(m_1x_tm_2)=\sigma_0(m_1)\sigma_t(x_t)\sigma_0(m_2)$), and $\sigma_{t+s}(x_t\odot y_s)=\sigma_t(x_t)\sigma_s(y_s)$. In the \nbd{W^*}case we say $\sigma^\odot$ is \hl{normal}, if every $\sigma_t$ is \nbd{\sigma}weak with respect to the natural \nbd{\sigma}weak topologies. A covariant representation is \hl{nondegenerate}, if $\sigma_t(F_t)G$ is total in $G$ for every $t$. It is \hl{isometric}, if $\sigma_t(x_t)^*\sigma_t(y_t)=\sigma_0(\AB{x_t,y_t})$. Note that an isometric bimodule map $\sigma_t$ is completely contractive, automatically. In the \nbd{W^*}case, $\sigma_t$ is normal, if and only if the restriction of $\sigma_0$ to the \nbd{W^*}subalgebra $\ol{F_M}^s$ of $M$ generated by the \hl{range ideal} $F_M:=\cls\AB{F,F}$ is normal. A \hl{(normal) isometric dilation} of a (normal) covariant representation $\sigma^\odot$ on $G$ is a (normal) covariant isometric representation $\tau^\odot$ on a Hilbert space $H\supset G$ such that $p_G\tau_t(\bullet)p_G=\sigma_t$, where $p_G\in\sB(H)$ is the projection onto $G$.

A \it{nondegenerate} covariant representation is what Muhly and Solel \cite{MuSo02} call \it{fully coisometric}. In \cite[Theorem 3.3]{MuSo98} they proved the result for a \it{single} correspondence, what in the notation of Theorem \ref{CRdilthm} is just the discrete case $\bS=\N_0$. (In this case $F^\odot=\bfam{F^{\odot n}}_{n\in\N_0}$ is just the product system generated by the single correspondence $F$.) In \cite[Theorem 3.7]{MuSo02} they proved the continuous time case $\bS=\R_+$ for a \it{nondegenerate} representation $\sigma^\odot$. (We should note that the word ``discrete'' in their theorem does not refer to the indexing semigroup $\bS$, but simply is to emphasize that there are no continuity or measurability conditions on the product system. We prefer to use the term \hl{algebraic} product system, as ``discrete'', in our opinion, fits much better to the semigroup $\bS$.) In Theorem \ref{CRdilthm} we treat now the general case.

The heart of our proof is the following theorem about existence of a (weak) dilation of a CP-semigroup with a pre-assigned product system. It has an obvious extension to \nbd{W^*}modules, which, in a sense, is dual to Theorem \ref{CRdilthm}.

\bthm\label{CPdilthm}
Let $E^\odot=\bfam{E_t}_{t\in\bS}$ be a product system of correspondences over a unital \nbd{C^*}al\-ge\-bra $\cB$ and let $\xi^\odot=\bfam{\xi_t}_{t\in\bS}$ be a contractive (that is, $\AB{\xi_t,\xi_t}\le\U$) unit for $E^\odot$. Then there exists a unique strict weak dilation $(E,\vt,\xi)$ of the CP-semigroup $T=\bfam{T_t}_{t\in\bS}$ $(T_t:=\AB{\xi_t,\bullet\xi_t})$ fulfilling:
\begin{enumerate}
\item
The product system associated with the strict \nbd{E}semigroup $\vt=\bfam{\vt_t}_{t\in\bS}$ on $\sB^a(E)$ is $E^\odot$.

\item\label{CPd2}
$E$ is generated by $E^\odot$ in the sense that $\bigcup_{t\in\bS}\vt_t(\xi\xi^*)E$ is total in $E$.
\end{enumerate}
\ethm

\noindent
Recall that a \hl{unit} for a $E^\odot$ is a family $\xi^\odot=\bfam{\xi_t}_{t\in\bS}$ of elements $\xi_t\in E_t$ with $\xi_0=\U$ that is multiplicative in the sense that $\xi_s\odot\xi_t=\xi_{s+t}$. Clearly, for every (contractive) unit, the mappings $T_t$ define a semigroup of (contractive) CP-maps on $\cB$. An \hl{\nbd{E}semigroup} is a semigroup of endomorphisms of a \nbd{*}algebra. If these endomorphisms are unital, then we say \hl{\nbd{E_0}semigroup}. An \nbd{E}semigroup on $\sB^a(E)$ is \hl{strict}, if every endomorphism $\vt_t$ is strict on bounded subsets. (Equivalently, the action $\vt_t(\sK(E))$ of the compacts $\sK(E)$ is already nondegenerate on $\vt_t(\U)E$.) Finally, a \hl{weak dilation} of a CP-semigroup $T$ of contractions $T_t$ on $\cB$ is a triple $(E,\vt,\xi)$ consisting of a Hilbert \nbd{\cB}module, a strict \nbd{E}semigroup $\vt=\bfam{\vt_t}_{t\in\bS}$ on $\sB^a(E)$ and a \hl{unit vector} $\xi\in E$ (that is, $\AB{\xi,\xi}=\U$) such that $T_t(b)=\AB{\xi,\vt_t(\xi b\xi^*)\xi}$. The term ``weak'' refers to that the representations $b\mapsto\vt_t(\xi b\xi^*)$ of $\cB$ (that is, the \hl{Markov flow} associated with the dilation) are, usually, nonunital; see Bhat and Skeide \cite{BhSk00}. We recall in Section \ref{dilSEC} what we mean by the product system associated with an \nbd{E}semigroup.

We prove Theorem \ref{CPdilthm} in Section \ref{dilSEC}, as stated, in a version for \nbd{C^*}correspondences over a unital \nbd{C^*}al\-ge\-bra. It extends easily to \nbd{W^*}correspondences. In Section \ref{comSEC} we recall what the \it{commutant} of a product system of \nbd{W^*}correspondences (or better, of von Neumann correspondences) is. We make some effort to point out that the correspondence between a product system and its commutant is bijective when applied to categories of \it{concrete} von Neumann correspondences. We establish a couple of bijections (Theorem \ref{EE'thm}) between structures present in a product system and structures present in its commutant system. In particular, we show that after passing to the \it{commutant system} of $E^\odot$ Theorem \ref{CPdilthm} translates into the \nbd{W^*}version of Theorem \ref{CRdilthm} in the special case when $\sigma_0$ is faithful and nondegenerate. In Section \ref{proofSEC} we explain how both the general \nbd{W^*}case and the \nbd{C^*}case boil down to Theorems \ref{CPdilthm} and \ref{EE'thm}.

\lf
We mention that Muhly and Solel (private communication) give a direct formulation of the construction of the isometric dilation in Theorem \ref{CRdilthm} without passing through commutants. However, their proof only reduces the problem to their dilation result for \it{nondegenerate} covariant representations \cite[Theorem 3.7]{MuSo02}. We should like to say that, assuming the basic and (easy to verify) facts about commutants of von Neumann correspondences, the proof we give here is self-contained.

It is very well possible to read the three following sections in reverse order, that is, first reducing the general statement to the case of faithful nondegenerate $\sigma_0$, then showing equivalence of the statement of Theorem \ref{CRdilthm} in this case with the statement of Theorem \ref{CPdilthm} and, finally, proving Theorem \ref{CPdilthm}. Of course, reading in this order makes less clear that Theorem \ref{CPdilthm} is a result that is independent of the remainder.

Note that, in Theorem \ref{CRdilthm}, as compared with our usual convention, we have switched the order of the letters $s,t$ for time arguments. This may appear to be a total triviality. However, in the transition between a product system and its commutant in Section \ref{comSEC} also time-orders do change; see Remark \ref{coswitchrem}. And it is Theorem \ref{CPdilthm} which is directly related to our usual conventions, while Theorem \ref{CRdilthm}, the dual of Theorem \ref{CPdilthm}, has opposite order.

\lf\noindent
\bf{Acknowledgements.~}
Apart from support PRIN 2005 financed by the Italian MIUR and by the University of Molise, we would like to thank the referee having pointed out a fatal typo in Section \ref{dilSEC} of the first version, for having filled in some bibliographic gaps, and for encouraging to be more explicit about the involved categories in Section \ref{comSEC}.


\section{Proof of Theorem \ref{CPdilthm}}\label{dilSEC}

There are several ways to associate with an \nbd{E}semigroup a product system. These possibilities are all equivalent in that they produce isomorphic product systems. We will have to give a more specific meaning to this in Section \ref{comSEC}. In this section we are interested only in the special case of the product system \hl{associated} with a strict \nbd{E}semigroup $\vt$ acting on $\sB^a(E)$ where $E$ is a Hilbert \nbd{\cB}module with a \hl{unit vector} $\xi\in E$ (that is, $\AB{\xi,\xi}=\U$). In this case the product system associated with an \nbd{E}semigroup (or better: one representative of the isomorphism class) may be obtained in exactly the same way as in Skeide \cite{Ske02}: Put $E_t:=\vt_t(\xi\xi^*)E$. Then $E_t$ becomes a correspondence over $\cB$ by defining the left action $bx_t:=\vt_t(\xi b\xi^*)x_t$ of $b\in\cB$. By $x\odot y_t\mapsto\vt_t(x\xi^*)y_t$ we define a unitary $E\odot E_t$ onto $\vt_t(\U)E\subset E$ giving back $\vt_t$ as $\vt_t(a)=a\odot\id_{E_t}\in\sB^a(\vt_t(\U)E)\subset\sB^a(E)$. Moreover, the restriction of this unitary to $E_s\odot E_t\subset E\odot E_t$ defines a bilinear unitary $u_{s,t}\colon E_s\odot E_t\rightarrow E_{s+t}$ and the family of all these bilinear unitaries defines a product system structure. If there is no danger of confusion we will suppress the $u_{s,t}$ and simply write $E_s\odot E_t=E_{s+t}$. Note that $E_t\subset E$, so that it makes sense to speak about $E$ being generated by $E^\odot$ in the sense of Theorem \ref{CPdilthm}\eqref{CPd2}. The generalization of \cite{Ske02} from \nbd{E_0}semigroups to \nbd{E}semigroups has been discussed by Bhat and Lindsay \cite{BhLi05}. We may also obtain a version for nonunital $\cB$ (that is, in particular, without unit vector) following the methods from Muhly, Skeide and Solel \cite{MSS06} as explained in the case of \nbd{E_0}semigroups in \cite{Ske04p}. However, note that in this case large parts of the proof of Theorem \ref{CPdilthm} will not work.

The idea to prove Theorem \ref{CPdilthm} is to mimic the construction in \cite[Section 8]{BhSk00} which treats the case when $E^\odot$ is the \hl{GNS-system} of the CP-semigroup $T$, that is, the unit $\xi^\odot$ that gives back $T$ as $\AB{\xi_t,\bullet\xi_t}$ generates the whole product system. The construction in \cite{BhSk00} is done by passing to the unitalized CP-semigroup $\wt{T}$ on the unitalization $\wt{\cB}=\C\wt{\U}\oplus\cB$. Then the so-called \it{minimal dilation} of $\wt{T}$ is constructed. This construction involves an inductive limit of correspondences (called the first inductive limit in \cite{BhSk00} or the two-sided inductive limit in \cite{BBLS04}) giving the product system and an inductive limit of right Hilbert modules (called the second inductive limit in \cite{BhSk00} or the one-sided inductive limit in \cite{BBLS04}) giving the right module on which the \nbd{E_0}semigroup dilating $\wt{T}$ lives. At the moment of the proof in \cite{BhSk00} these constructions are already known and simply used. The major part of the work in \cite[Section 8]{BhSk00} consists in identifying how the dilation of $T$ itself sits inside the dilation of $\wt{T}$.

Here we have to redo the first inductive limit (or better to replace it with something similar) to obtain a product system of correspondences over $\wt{\cB}$ with a unital unit that allows, then, to do the second inductive limit, also yielding a dilation (just not necessarily the minimal one) of $\wt{T}$. The work to show how this dilation of $\wt{T}$ contains a dilation of $T$ is very similar to \cite{BhSk00}. We apologize for that with respect to \cite{BhSk00} we find it convenient to change the order of certain components in column vectors. (For instance, since $\cB$ is unital, $\wt{\cB}$ is isomorphic to the \nbd{C^*}algebraic direct sum $\C\oplus\cB$. Here we shall write the component $\C$ as upper component, while in \cite{BhSk00} we wrote it as lower component.) We also mention a typo in \cite[Theorem 8.4]{BhSk00} where we wrote accidentally that $\vt$ is an \nbd{E_0}semigroup. (Of course, if $T$ is nonunital, then $\vt$ constructed in \cite{BhSk00} is definitely not an \nbd{E_0}semigroup but only an \nbd{E}semigroup.) We also should like to say that the version in \cite{BhSk00} is formulated for pre-Hilbert modules and extends easily to Hilbert modules (and, likewise, to \nbd{W^*} or von Neumann modules), while here we write immediately for Hilbert modules. This is mere convenience, and whatever we write down in this section would work also in the algebraic context without any completion.

A crucial point in \cite{BhSk00} was to identify the GNS-module of $\wt{T}_t$ and its cyclic vector $\wt{\xi}_t$ in terms of the GNS-module $\sE_t$ and cyclic vector $\xi_t$ of $T_t$.\footnote{The \hl{GNS-construction} for a CP-map $T$ on (or between) unital (pre-)\nbd{C^*}algebra(s) is due to Paschke \cite[Theorem 5.2]{Pas73}. The result is a correspondence $\sE$ and vector $\xi\in\sE$ such that $T$ is recovered as $T=\AB{\xi,\bullet\xi}$. If $\xi$ is cyclic for $\sE$, that is if $\xi$ generates $\sE$ as a correspondence, then we speak of the \hl{GNS-module}, as in this case everything is determined up to suitable unitary equivalence.} Here we do the same for $E_t$ and $\xi_t$, just that now $E_t$ may be bigger than $\sE_t$. So put $\wh{\xi}_t:=\sqrt{\,\wt{\U}-\AB{\xi_t,\xi_t}}\in\wt{\cB}$ and denote by $\wh{E}_t:=\ol{\wh{\xi}_t\wt{\cB}}$ the closed right ideal in $\wt{\cB}$ generated by $\wh{\xi}_t$. Turn the Hilbert \nbd{\wt{\cB}}module $\wh{E}_t$ into a correspondence over $\wt{\cB}$ by putting $b\wh{\xi}_t=0$ for every $b\in\cB$ (and, of course, $\wt{\U}~\wh{\xi}_t=\wh{\xi}_t$). Moreover, $E_t$ is a correspondence over $\wt{\cB}$ in the unique (because nondegenerate) way. Define the \nbd{\wt{\cB}}correspondence $\wt{\sE}_t:=\wh{E}_t\oplus E_t$. Observe that $\wt{\xi}_t:=\wh{\xi}_t\oplus\xi_t$ is a unit vector in $\wt{\sE}_t$ and that $\AB{\wt{\xi}_t,\bullet\wt{\xi}_t}$ is nothing but the unitalization of $T_t$. (In fact, if $E_t$ was the GNS-module, then $\wt{\sE}_t$ would be the GNS-module of $\wt{T}_t$; see \cite{BhSk00}.)

As $E_s\odot\wh{E}_t=\zero$, we find
\beqn{
\wt{\sE}_s\odot\wt{\sE}_t
~=~
(\wh{E}_s\odot\wh{E}_t)\,\oplus\,(\wh{E}_s\odot E_t)\,\oplus\,(E_s\odot E_t),
}\eeqn
and, similarly,
\beq{\label{wtEtdef}
\wt{\sE}_{t_n}\odot\ldots\odot\wt{\sE}_{t_1}
~=~
\bigoplus_{k=0}^n\wh{E}_{t_n}\odot\ldots\odot\wh{E}_{t_{k+1}}\odot E_{t_k}\odot\ldots\odot E_{t_1}.
}\eeq
Accordingly,
\beqn{
\wt{\xi}_{t_n}\odot\ldots\odot\wt{\xi}_{t_1}
~=~
\bigoplus_{k=0}^n\wh{\xi}_{t_n}\odot\ldots\odot\wh{\xi}_{t_{k+1}}\odot\xi_{t_k}\odot\ldots\odot\xi_{t_1}.
}\eeqn
Of course, this is the tensor product of unit vectors and, therefore, itself a unit vector. (A direct verification of this trivial fact alone from the right-hand side would be quite tedious.) Rather than the lattices of interval partitions of $\LO{0,t}$, we use the lattices
\beqn{
\bJ_t
~:=~
\BCB{\,\et=(t_n,\ldots,t_1)\colon ~n\in\N\,,\,t_i>0\,,\,t_n+\ldots+t_1=t\,}
}\eeqn
$(t>0)$ as in \cite{BhSk00}. We define a partial order on $\bJ_t$ by $\es\le\et$, if there are $\es_j\in\bJ_{s_j}$, $\es=(s_m,\ldots,s_1)\in\bJ_t$ such that $\et=\es_m\smallsmile\ldots\smallsmile\es_1$, where the \hl{join} $\smallsmile$ of two tuples $\es=(s_m,\ldots,s_1)\in\bJ_s,\et=(t_n,\ldots,t_1)\in\bJ_t$ is defined as $\es\smallsmile\et:=(s_m,\ldots,s_1,t_n,\ldots,t_1)\in\bJ_{s+t}$. See \cite[Observation 4.2]{BhSk00} for why we find this lattice more useful than the lattice of interval partitions (to which it is isomorphic). For $\et=(t_n,\ldots,t_1)\in\bJ_t$ it follows that
\baln{
\beta_{\et(t)}
&
\colon
&
\wh{\xi}_t
&
~\longmapsto~
\bigoplus_{k=0}^{n-1}\wh{\xi}_{t_n}\odot\ldots\odot\wh{\xi}_{t_{k+1}}\odot\xi_{t_k}\odot\ldots\odot\xi_{t_1}
&
&
\text{and}
&
E_t
&
~\xrightarrow{~~\cong~~}~
E_{t_n}\odot\ldots\odot E_{t_1}
}\ealn
defines a bilinear isometric embedding $\wt{\sE}_t\rightarrow\wt{\sE}_{t_n}\odot\ldots\odot\wt{\sE}_{t_1}=:\wt{\sE}_\et$ sending $\wt{\xi}_t$ to $\wt{\xi}_{t_n}\odot\ldots\odot\wt{\xi}_{t_1}$. (To check that this mapping is isometric on $\wh{\xi}_t$ simply observe that the missing term $k=n$ in the sum has ``square length'' $\AB{\xi_t,\xi_t}$, so that the sum from $k=0$ to $k=n-1$ has ``square length'' $\wt{\U}-\AB{\xi_t,\xi_t}=\AB{\wh{\xi}_t,\wh{\xi}_t}$ as it should, and that the left action of $b\in\cB$ on $\beta_{\et(t)}(\wh{\xi}_t)$ gives $0$ as it should.) For $\et\ge\es$ (so that $\et=\es_m\smallsmile\ldots\smallsmile\es_1$, $\es_j\in\bJ_{s_j}$, $\es=(s_m,\ldots,s_1)\in\bJ_t$) we put $\beta_{\et\es}:=\beta_{\es_m(s_m)}\odot\ldots\odot\beta_{\es_1(s_1)}$. As in \cite[Section 4]{BhSk00} the $\beta_{\et\es}\colon\wt{\sE}_\es\rightarrow\wt{\sE}_\et$ form an inductive system, the inductive limits $\wt{E}_t:=\limind_{\et\in\bJ_t}\wt{\sE}_\et$ form a product system $\wt{E}^\odot$ and the $\wt{\xi}_t\in\wt{\sE}_t=\wt{\sE}_{(t)}\subset\wt{E}_t$ form a unital unit $\wt{\xi}^\odot$ for $\wt{E}^\odot$. (Of course, $\wt{E}_s\odot\wt{E}_t$ imbeds into $\wt{E}_{s+t}$, and surjectivity follows from the fact that every tuple in $\bJ_{s+t}$, by adding a single time point if necessary, has a refinement of the form $\es\smallsmile\et$ with $\es\in\bJ_s,\et\in\bJ_t$. In addition to \cite[Section 4]{BhSk00}, see also Barreto, Bhat, Liebscher and Skeide \cite[Section 4.3]{BBLS04} and Skeide \cite{Ske03c,Ske06d} for similar two-sided limits. See \cite[Appendix A]{BhSk00} for details about inductive limits of Hilbert modules and correspondences.)

Clearly, $\cls\wt{\cB}\,\wt{\xi}_t\wt{\cB}$ contains $\wh{E}_t=(\wt{\U}-\U)\wt{\sE}_t$. This shows that $\wt{E}^\odot$, as a product system, is generated by $E^\odot$ and $\wt{\xi}^\odot$ (in the sense that there is no proper subsystem of $\wt{E}^\odot$ containing $E^\odot$ and $\wt{\xi}^\odot$). Also, since $\wh{\xi}_t=(\wt{\U}-\U)\wt{\xi}_t$, we find that
\beqn{
\wh{\xi}_{s+t}
~=~
(\wt{\U}-\U)\wt{\xi}_s\odot\wt{\xi}_t
~=~
\wh{\xi}_s\odot\wt{\xi}_t
~=~
\wh{\xi}_s\odot\wh{\xi}_t+\wh{\xi}_s\odot\xi_t
}\eeqn
so that everything in $\wt{E}_t$ that lies in the complement of $E_t$ lies in the span of elements of the form $\wh{\xi}_{t-s}\odot x_s$ $(0\le s<t,x_s\in E_s)$. Note that by \eqref{wtEtdef}, $\U\wt{E}_t$ is just $E_t$.

From the unital unit $\wt{\xi}^\odot$ we construct an inductive limit $\wt{E}=\limind_{t}\wt{E}_t$ as in \cite[Section 5]{BhSk00} with the help of the isometric embeddings $\wt{E}_s\rightarrow\wt{\xi}_s\odot\wt{E}_t\subset\wt{E}_{s+t}$. (In \cite{BhSk00} we discussed the minimal case, but the construction works without changing a word also in the general case. See, for instance, \cite[Section 4.4]{BBLS04} or Skeide \cite{Ske03b}.) We have the factorization $\wt{E}=\wt{E}\odot\wt{E}_t$ such that $\wt{\vt}(a)=a\odot\id_{\wt{E}_t}$ defines an \nbd{E_0}semigroup acting on $\sB^a(\wt{E})$, having the product system $\wt{E}^\odot$. Further, we have a unit vector $\wt{\xi}=\wt{\xi}_t\in\wt{E}_t\subset\wt{E}$, satisfying $\wt{\xi}=\wt{\xi}\odot\wt{\xi}_t$, so that $(\wt{E},\wt{\vt},\wt{\xi})$ is a dilation of $\wt{T}$.

Next we show how the dilation of $\wt{T}$ restricts to a dilation of $T$. We proceed as in \cite[Section 8]{BhSk00}. We define $E:=\wt{E}\U$ and observe that $E$ is a Hilbert \nbd{\cB}module with a unit vector $\xi:=\wt{\xi}\U$. As multiplication with $\U$ from the right defines a central projection $p$ in $\sB^a(\wt{E})$ onto $E$, an operator $a\in\sB^a(\wt{E})$ is in $\sB^a(E)$, if and only if $pa(=ap)=a$. So, for $a\in\sB^a(E)$ it follows that
\beqn{
\wt{\vt}_t(a)
~=~
a\odot\id_{\wt{E}_t}
~=~
(pa)\odot\id_{\wt{E}_t}
~=~
a\odot\id_{E_t},
}\eeqn
where, by slight abuse of notation, we denote the projection onto $E_t$ in $\sB^a(\wt{E}_t)$ (that is, left multiplication with $\U\in\sB$) by $\id_{E_t}$. This shows that $\wt{\vt}_t$ leaves $\sB^a(E)$ invariant. Denote by $\vt$ the restriction of $\wt{\vt}$ to $\sB^a(E)$. Then,
\beqn{
\vt_t(\xi\xi^*)E
~=~
\wt{\vt}_t(\wt{\xi}\U\wt{\xi}^*)\wt{E}\U
~=~
\wt{\vt}_t(\wt{\xi}\U\wt{\xi}^*)\wt{\vt}_t(\wt{\xi}\wt{\xi}^*)\wt{E}\U
~=~
\U\wt{E}_t\U
~=~
E_t\U
~=~
E_t,
}\eeqn
so that the product system of $\vt$ is $E^\odot$. (It is an easy exercise to show that the identification $\wt{E}_s\odot\wt{E}_t=\wt{E}_{s+t}$ restricted to elementary tensors from $E_s\odot E_t$ gives the correct identification $E_s\odot E_t=E_{s+t}$.) Also $\wt{\vt}_t(\wt{\xi}b\wt{\xi}^*)=\wt{\vt}_t(\xi b\xi^*)=\vt_t(\xi b\xi^*)$ and $\xi=\wt{\xi}\U=\wt{\xi}\odot\wt{\xi}_t\U$. Therefore, $\vt_t(\xi\xi^*)\xi=\xi\odot\xi_t=\xi_t$. That is, we obtain back the unit $\xi^\odot$ and $(E,\vt,\xi)$ is a dilation of $T$.

For showing Condition \ref{CPd2} in Theorem \ref{CPdilthm}, we first show that a dilation fulfilling Condition \ref{CPd2} is determined uniquely. Let consider an inner product (in $E$) of elements $x_t=\vt_t(\xi\xi^*)x\in E_t$ $(x\in E)$ and $y_s\odot z_t=\vt_s(\xi\xi^*)y\odot\vt_t(\xi\xi^*)z\in E_{s+t}$ $(y,z\in E)$. We find
\bmun{
\AB{x_t,y_s\odot z_t}
~=~
\BAB{\vt_t(\xi\xi^*)x\,,\,\vt_t\bfam{\vt_s(\xi\xi^*)y\xi^*}\vt_t(\xi\xi^*)z}
\\
~=~
\BAB{\vt_t\bfam{\vt_s(\xi\xi^*)\xi\xi^*}\vt_t(\xi\xi^*)x\,,\,\vt_t\bfam{\vt_s(\xi\xi^*)y\xi^*}\vt_t(\xi\xi^*)z}
~=~
\AB{\xi_s\odot x_t,y_s\odot z_t},
}\emun
so that these inner products (and, therefore, all inner products of $E$) can be calculated by using the product system structure of $E^\odot$ and the unit $\xi^\odot$.

It remains to show that the dilation we constructed fulfills Condition \ref{CPd2}. But this follows from totality of the elements $\wh{\xi}_{t-s}\odot x_s$ $(0\le s\le t,x_s\in E_s)$ in $\wt{E}_t$ and from $\wh{\xi}_t\U=\wt{\xi}_t\U-\U\wt{\xi}_t=\wt{\xi}_t\U-\U\wt{\xi}_t\U$ (see the crucial \cite[Observation 8.1]{BhSk00}). Indeed, the elements $\wt{\xi}\odot\wt{x}_t$ $(t\in\bS,\wt{x}_t\in\wt{E}_t)$ are dense in $\wt{E}$ so that the elements $\wt{\xi}\odot\wt{x}_t\U$ are dense in $E$. Now for all $0\le s\le t$ we have
\beqn{
\wt{\xi}\odot\wh{\xi}_{t-s}\odot x_s\U
~=~
\wt{\xi}\odot\wh{\xi}_{t-s}\U\odot x_s
~=~
\wt{\xi}\odot(\wt{\xi}_{t-s}\U-\U\wt{\xi}_{t-s})\odot x_s
~=~
\wt{\xi}\odot x_s-\wt{\xi}\odot\xi_{t-s}\odot x_s.
}\eeqn
Now $\wt{\xi}\odot x_s$ is in $E_s=\vt_s(\xi\xi^*)E$ and $\wt{\xi}\odot\xi_{t-s}\odot x_s$ is in $E_t=\vt_t(\xi\xi^*)E$. This shows that Condition \ref{CPd2} is fulfilled.

\brem
The discussion in \cite[Section 12]{BhSk00} shows how to adapt the arguments to von Neumann modules using the appendices of \cite{BhSk00}. Taking into account that every \nbd{W^*}module may be considered as a von Neumann module (by choosing a concrete representation of the underlying \nbd{W^*}algebra), the result holds also for \nbd{W^*}modules.
\erem

\brem
We would like to note that, like in \cite[Section 8]{BhSk00}, the \nbd{\C}linear codimension of $E$ in $\wt{E}$ is $1$. More precisely, $\Om:=\wt{\xi}(\wt{\U}-\U)$ is a vector with ``length'' $\AB{\Om,\Om}=\wt{\U}-\U$ such that $\wt{E}=\C\Om\oplus E$. This follows by looking at \eqref{wtEtdef} and from the fact that $\Om_t:=\wt{\xi}_t(\wt{\U}-\U)$ $(t>0)$ is a vector with ``length'' $\AB{\Om_t,\Om_t}=\wt{\U}-\U$ such that $\wt{\sE}_t=\C\Om_t\oplus\wt{\sE}_t\U$ and further $\wt{E}_t=\C\Om_t\oplus\wt{E}_t\U$. Finally, also $\Om^\odot=\bfam{\Om_t}_{t\in\bS}$ (with $\Om_0:=\wt{\U}$) is a unit for $\wt{E}^\odot$.
\erem

\section{Duality between dilations}\label{comSEC}

Von Neumann correspondences from $\cA$ to $\cB$ are in duality with von Neumann correspondences from $\cB'$ to $\cA'$ via the \it{commutant}, a functor that naturally generalizes the functor that sends a von Neumann algebra to its commutant. Also a product system of von Neumann \nbd{\cB}correspondences $E^\odot$ has a \it{commutant} $E'^\odot$, a product system of von Neumann \nbd{\cB'}correspondences. Under this duality contractive units $\xi^\odot$ for $E^\odot$ correspond to covariant representations $\sigma'^\odot$ of $E'^\odot$ with faithful $\sigma'_0$ and the contractive CP-semigroup $T_t=\AB{\xi_t,\bullet\xi_t}$ may, likewise, be reconstructed from $\sigma'^\odot$, simply by going back to the (unique) unit $\xi^\odot$ on the side of $E^\odot$. Also, \nbd{E}semigroups with associated product system $E^\odot$, under commutant, correspond to isometric representations of the commutant system $E'^\odot$. See Theorem \ref{EE'thm}(\ref{EE'1}-\ref{EE'2}).

The scope of this section is to combine these dualities in order to establish a duality between, on the one hand, weak dilations with pre-assigned product system $(E^\odot,\xi^\odot)$ of a normal CP-semigroup $T$ determined by a unit $\xi^\odot$ in the product system $E^\odot$ of von Neumann correspondences over a von Neumann algebra $\cB\subset\sB(G)$ and, on the other hand, isometric dilations of the covariant representations $\sigma'^\odot$ of the commutant system $E'^\odot$ on $G$ with $\sigma'_0=\id_{\cB'}$ that is associated with that CP-semigroup. See Theorem \ref{EE'thm}\eqref{EE'4}. This duality translates, then, the existence result Theorem \ref{CPdilthm} into the existence result Theorem \ref{CRdilthm} in the special case $\sigma'_0=\id_{\cB'}$.

The correspondence between a von Neumann algebra and its commutant is bijective. In order that this desirable property remains true for commutants of von Neumann correspondences (not degenerating to an equivalence), we have to choose our categories carefully. The correct category that allows to view the commutant as a bijective functor is the category of \it{concrete} von Neumann correspondences; Skeide \cite{Ske06b}. In the sequel, we discuss only the case relevant to us, namely, correspondences over $\cB$. In order that all correspondences stated in Theorem \ref{EE'thm} are true one-to-one correspondences (and not just up to isomorphism or equivalence) we will have to come back also to the problem mentioned in the beginning of Section \ref{dilSEC}, namely, to indicate more specifically when we consider a product system as the one associated with an \nbd{E}semigroup.

\lf
Before we can speak about concrete von Neumann correspondences, we have to speak about concrete von Neumann modules. Recall that a von Neumann algebra is a strongly closed \nbd{*}algebra $\cB\subset\sB(G)$ of operators acting nondegenerately on a Hilbert space $G$. As usual, by $\cB'\subset\sB(G)$ we denote the commutant of $\cB$. Similarly, a \hl{concrete von Neumann \nbd{\cB}module} is a subset $E$ of $\sB(G,H)$, where $H$ is another Hilbert space, such that
\begin{enumerate}
\item\label{cvnM1}
$E$ is a right \nbd{\cB}submodule of $\sB(G,H)$, that is, $E\cB\subset E$,

\item\label{cvnM2}
$E$ is a pre-Hilbert \nbd{\cB}module with inner product $\AB{x,y}=x^*y$, that is, $E^*E\subset\cB$,

\item\label{cvnM3}
$E$ acts nondegenerately on $G$, that is, $\cls EG=H$, and

\item\label{cvnM4}
$E$ is strongly closed.
\end{enumerate}
If we wish to underline the Hilbert space $H$, we will also write the pair $(E,H)$ for the concrete von Neumann \nbd{\cB}module.  One may show (see Skeide \cite{Ske00b,Ske05c}) that a subset $E$ of $\sB(G,H)$ fulfilling \ref{cvnM1}--\ref{cvnM3} (that is, $E$ is a concrete pre-Hilbert \nbd{\cB}module) is a concrete von Neumann \nbd{\cB}module, if and only if $E$ is self-dual, that is, if and only if $E$ is a \nbd{W^*}module over the von Neumann algebra $\cB\subset\sB(G)$ considered as a \nbd{W^*}algebra. By $\cvN_\cB$ we denote the \hl{category of concrete von Neumann \nbd{\cB}modules} with the adjointable (in the usual sense) and, therefore, bounded maps $a\in\sB^a(E_1,E_2)$ as morphisms. The definition of concrete von Neumann modules and their category is due to Skeide \cite{Ske06b}.

Identifying $xg\in H$ with $x\odot g\in E\odot G$, we see from \ref{cvnM3} that $H$ and $E\odot G$ are canonically isomorphic.

\brem
In fact, if $E$ is a pre-Hilbert module over a pre-\nbd{C^*}algebra $\cB\subset\sB(G)$, then one may construct the Hilbert space $E\odot G$ with an embedding $x\mapsto L_x\in\sB(G,E\odot G)$ where we put $L_xg:=x\odot g$, transforming $E$ into a concrete pre-Hilbert \nbd{\cB}module $(E,E\odot G)$. For a von Neumann algebra $\cB\subset\sB(G)$ we defined in Skeide \cite{Ske00b} that $E$ is a von Neumann \nbd{\cB}module, if its image in $\sB(G,E\odot G)$ is strongly closed. Of course, in that way also a \nbd{W^*}module over a \nbd{W^*}algebra $M$ may be turned into a von Neumann module after choosing a faithful normal unital representation of $M$ on a Hilbert space $G$, thus, turning $M$ into a von Neumann algebra.
\erem

Giving $E$ as a subset of $\sB(G,H)$ from the beginning, is crucial for that the commutant, later on, will be bijective. However, the fact that $H$ is canonically isomorphic to the tensor product $E\odot G$ is by far more inspiring from the algebraic point of view. For instance, by well-known standard results certain operators on the factors $E$ or $G$ of a tensor product embed into the operators on the tensor product $E\odot G$ via amplification, while the corresponding action on $H$ is much less intuitive. So, every adjointable operator $a\in\sB^a(E_1,E_2)$ amplifies to an operator $a\odot\id_G\in\sB(E_1\odot G,E_2\odot G)$. Consequently, $a$ gives rise to and is determined uniquely by an operator in $\sB(H_1,H_2)$ that acts as $x_1g\mapsto(ax_1)g$. We shall denote this operator by the same symbol $a$ and indentify in that way $\sB^a(E_1,E_2)$ as a subset of $\sB(H_1,H_2)$. It is easy to show that $\sB^a(E_1,E_2)$ is strongly closed in $\sB(H_1,H_2)$. In particular, $\sB^a(E)\subset\sB(H)$ is a von Neumann algebra acting on $H$.

The operators on the second factor in $E\odot G$ that embed into $\sB(E\odot G)$ are the bilinear, that is, the \nbd{\cB}\nbd{\C}linear operators on $G$. Of course, $\sB^{bil}(G)=\cB'$ is nothing but the commutant of $\cB$. So, the (clearly, normal and nondegenerate) representation $b'\mapsto\id_E\odot b'$ of $\cB'$ on $E\odot G$ gives rise to a normal nondegenerates representation $\rho'$ of $\cB'$ on $H$ which acts as $\rho'(b')xg=xb'g$. We call $\rho'$ the \hl{commutant lifting} associated with $E$.

\brem
In both cases, the action of $a\in\sB^a(E_1,E_2)$ as an operator in $\sB(H_1,H_2)$ and the action of $\rho'(b')$ on $H$, there is no problem in showing that these operators are well-defined. But only the tensor product pictures $a\odot\id_G$ and $\id_E\odot b'$ explain where the operators come from and why it is natural to write them down.
\erem

Let us return to the concrete von Neumann \nbd{\cB}module $(E,H)$. From the commutant lifting $\rho'$ we obtain back $E$ as the space
\beq{\label{CB'def}
E
~=~
C_{\cB'}(\sB(G,H))
~:=~
\bCB{x\in\sB(G,H)\colon\rho'(b')x=xb'~(b'\in\cB')}
}\eeq
of \it{intertwiners} for the natural actions of $\cB'$. (This was known already to Rieffel \cite{Rie74a}. See \cite{Ske05c} for a proof by simply calculating the double commutant of the linking von Neumann algebra in $\sB(G\oplus H)$. This proof also shows that the commutant $\rho'(\cB')'$ of the range of $\rho'$ in $\sB(H)$ may be identified with the von Neumann algebra $\sB^a(E)\subset\sB(H)$. By doing the computation for $E=E_1\oplus E_2$ one also shows that $\sB^a(E_1,E_2)$ is just $\sB^{bil}(H_1,H_2)$, the space of operators that \it{intertwine} the commutant liftings $\rho'_2$ and $\rho'_1$.) Conversely, if $(\rho',H)$ is a normal nondegenerate representation of $\cB'$ on the Hilbert space $H$, then $E:=C_{\cB'}(\sB(G,H))$ as in \eqref{CB'def} defines a concrete von Neumann \nbd{\cB}module in $\sB(G,H)$, which gives back $\rho'$ as commutant lifting. (The only critical task, nondegeneracy in Condition \ref{cvnM3}, follows from Muhly and Solel \cite[Lemma 2.10]{MuSo02}.) We find that
\bal{\label{MRfun}
(E,H)
&~\longleftrightarrow~
(\rho',H)
&
a\in\sB^a(E_1,E_2)
&~\longleftrightarrow~
a\in\sB^{bil}(H_1,H_2)
}\eal
establishes a bijective functor between the category $\cvN_\cB$ of concrete von Neumann \nbd{\cB}modules and the \hl{category $_{\cB'}\cvN$ of normal nondegenerate representations of $\cB'$} with the intertwiners $\sB^{bil}(H_1,H_2)$ as morphisms. (The preceding correspondence was established in Skeide \cite{Ske03c} as an equivalence between the category von Neumann \nbd{\cB}modules and $_{\cB'}\cvN$. As a von Neumann \nbd{\cB}module $E$, first, must be turned into a concrete von Neumann \nbd{\cB}modules $(E,E\odot G)$, the correspondence is not bijective but only an equivalence. The precise formulation above, where the functor is, really, bijective and not only an equivalence, is due to \cite{Ske06b}.)

\lf
A \hl{concrete von Neumann correspondence} over a von Neumann algebra $\cB$ is a concrete von Neumann \nbd{\cB}module $(E,H)$ with a left action of $\cB$ such that $\rho\colon\cB\rightarrow\sB^a(E)\rightarrow\sB(H)$ defines a normal (nondegenerate, of course) representation of $\cB$ on $H$. We call $\rho$ the \hl{Stinespring representation} associated with $E$.

\brem
If $E$ is the GNS-module of a (normal) CP-map, then $\rho$ is, indeed, the Stinespring representation, while $\rho'$ is (a restriction of) the representation constructed by Arveson \cite{Arv69} in the section called ``lifting commutants''.
\erem

By $_\cB\cvN_\cB$ we denote the \hl{category of concrete von Neumann correspondences from $\cB$ to $\cB$} with the the bilinear adjointable maps $a\in\sB^{a,bil}(E_1,E_2)$ as morphisms. (For adjointable maps, only left \nbd{\cB}linearity has to be checked.) We observe that $\rho(\cB)\subset\sB^a(E)=\rho'(\cB')'$, that is, $\rho'$ and $\rho$ have mutually commuting ranges. As this is very close to correspondences in the sense of Connes \cite{Con80p} (if $\cB$ is in \it{standard form}, then $\cB'\cong\cB^{op}$), we introduce the \hl{category of concrete Connes correspondences} $_\cB\ecC_\cB$ whose objects are triples $(\rho',\rho,H)$ such that $\rho'$ and $\rho$ are a pair of normal nondegenerate representations of $\cB'$ and of $\cB$, respectively, on $H$ with mutually commuting ranges, and with those maps in $\sB(H_1,H_2)$ as morphisms that intertwine both actions that of $\cB'$ and that of $\cB$. Extending the correspondence between concrete von Neumann \nbd{\cB}modules and representations of $\cB'$, we find a find bijective functor between the category of concrete von Neumann \nbd{\cB}correspondences $(E,H)$ and the category of concrete Connes correspondences $(\rho',\rho,H)$. In \cite{Ske03c} we observed this as an equivalence for von Neumann correspondences, while the bijective version for concrete von Neumann correspondences is from \cite{Ske06b}.

A last almost trivial observation (once again in \cite{Ske03c} up to equivalence and in \cite{Ske06b}, really, bijective) consists in noting that in the representation picture the roles of the representations $\rho'$ and $\rho$ are absolutely symmetric. That is, $_\cB\ecC_\cB\cong{_{\cB'}}\ecC_{\cB'}$. Therefore, if we switch $\cB$ and $\cB'$, that is, if we interprete $\rho$ as commutant lifting of $\cB$, the commutant of $\cB'$, and $\rho'$ as Stinespring representation of $\cB'$, by
\beq{\label{CBdef}
E'
~:=~
C_\cB(\sB(G,H))
~:=~
\bCB{x'\in\sB(G,H)\colon\rho(b)x'=x'b~(b\in\cB)}
}\eeq
we obtain a von Neumann \nbd{\cB'}module which is turned into a von Neumann \nbd{\cB'}correspondence by defining a left action via $\rho'$. We call $E'$ the \hl{commutant} of $E$. The commutant is a bijective functor from the category of concrete von Neumann \nbd{\cB}correspondences onto the category of concrete von Neumann \nbd{\cB'}correspondences (in each case with the bilinear adjointable maps as morphisms that are, really, the same algebra $\sB^a(E)\cap\sB^a(E')=\rho'(\cB')'\cap\rho(\cB)'$ of operators in $\sB(H)$). Obviously, $E'':=(E')'=E$.

\brem\label{gencommrem}
Muhly and Solel \cite{MuSo04} have discussed a version of the commutant for \nbd{W^*}al\-ge\-bras, called \it{\nbd{\sigma}dual}, where $\sigma$ is a faithful representation of the underlying \nbd{W^*}algebra, that must be chosen, and the \nbd{\sigma}dual depends on $\sigma$ (up to Morita equivalence of correspondences \cite{MuSo05}). The extension to correspondences from $\cA$ to $\cB$ was first done in the setting of \nbd{\sigma}duals in \cite{MuSo05}. In \cite{Ske06b} we discussed the version for von Neumann algebras and (concrete) von Neumann correspondences.

We remark that the functor $\cvN_\cB\leftrightarrow{_{\cB'}}\cvN$ in \eqref{MRfun} fits canonically into the setting of the commutant functor as $_\C\cvN_\cB\longleftrightarrow{_{\cB'}}\cvN_{\C'}$, if we consider $\C=\C'\subset\sB(\C)=\C$ as a von Neumann algebra.
\erem

The tensor product of Connes correspondences is tricky to describe in terms that do not explicitly involve the von Neumann correspondences to which they correspond. It requires that the von Neumann algebra is a \nbd{W^*}algebra in standard form and parts from Tomita-Takesaki theory and the result depends manifestly on the choice of a normal semifinite weight; see, for instance, Takesaki \cite[Section IX.3]{Tak03a}. Also the tensor product of \nbd{W^*}correspondences, although definitely less involved, still has the problem that the usual tensor product must be completed in a suitable \it{\nbd{\sigma}topology}, and this topology is defined rather \it{ad hoc}.

The tensor product two of von Neumann correspondences $E_1$ and $E_2$ is easy to obtain (and canonical up to unitary equivalence): Simply construct $E_1\odot E_2\odot G$ and determine the strong closure of $E_1\,\ul{\odot}\;E_2$ in $\sB(G,E_1\odot E_2\odot G)$ or, equivalently, determine the intertwiner space $C_{\cB'}(\sB(G,E_1\odot E_2\odot G))$, a purely algebraic problem, like determining the double commutant of a \nbd{*}algebra of operators. Up to canonical isomorphism it is not important whether we construct first $E_1\,\ul{\odot}\;E_2$ and then $(E_1\,\ul{\odot}\;E_2)\odot G$ or first $E_2\odot G$ and then $E_1\odot(E_2\odot G)$. If we have concrete von Neumann correspondences $(E_1,H_1)$ and $(E_2,H_2)$ it occurs to be more adapted to construct $E_1\odot H_2$ as the space $H_2$, canonically isomorphic to $E_2\odot G$, is given from the beginning. By slight abuse of notation we shall denote the concrete von Neumann correspondence obtained in that way by $E_1\odot E_2\subset\sB(G,E_1\odot H_2)$, using the same symbol $\odot$ as for the tensor product of \nbd{C^*}correspondences. Anyway, no matter how we obtained $E_1\odot E_2\odot G$, as $(E_1\,\ul{\odot}\;E_2)\odot G$, as $E_1\odot(E_2\odot G)$ or as $E_1\odot H_2$, to fix an isomorphism from the concrete von Neumann correspondence $(E_1\odot E_2,E_1\odot E_2\odot G)$ to a concrete von Neumann correspondence $(F,K)$ simply means to fix a unitary $u\in\sB(E_1\odot E_2\odot G,K)$ that intertwines both the commutant liftings of $\cB'$ and the Stinespring representations of $\cB$.

\lf
The notations established so far allow to state and prove Theorem \ref{EE'thm}\eqref{EE'1}, namely, that the commutant establishes a bijective functor between the category $\cvN^\odot_\cB$ of product systems of concrete von Neumann \nbd{\cB}cor\-re\-spond\-ences and the category $\cvN^\odot_{\cB'}$ of product systems of concrete von Neumann \nbd{\cB'}correspondences. A morphism between two objects $E^\odot$ and $F^\odot$ in $\cvN^\odot_\cB$ is a family $a^\odot=\bfam{a_t}_{t\in\bS}$ of maps $a_t\in\sB^{a,bil}(E_t,F_t)$ that fulfills $a_s\odot a_t=a_{s+t}$ and $a_0=\id_\cB$. Also Theorem \ref{EE'thm}\eqref{EE'3}, namely, that contractive units for $E^\odot$ correspond one-to-one with covariant normal representations of the commutant $E'^\odot$, may be stated and proved. (Actually, for this part it is not even necessary to speak about concrete von Neumann correspondences. It is true as soon as we fix a pair $E^\odot,E'^\odot$ of product systems of von Neumann correspondences that are commutants of each other up to isomorphism.)

So far, we associated with an \nbd{E}semigroup on $\sB^a(E)$ a product system, in the case that $E$ has a unit vector. But there are other ways to do this (see \cite{Ske03c,Ske04p}) that lead to canonically isomorphic, but definitely not equal product systems. In order that also the correspondence between \nbd{E}semigroups having a certain product system and isometric representations of the commutant of that product system becomes one-to-one, we have to indicate carefully what we understand by the opposite direction, that is, what it means that an \nbd{E}semigroup is associated with a product system. We shall say an \nbd{E}semigroup $\vt$ on $\sB^a(E)$ is \hl{associated} with a given product system $E^\odot$, if $\vt$ can be recovered as $\vt_t(a)=u_t(a\odot\id_t)u_t^*$ for a family of isometries $u_t\colon E\odot E_t\rightarrow E$ fulfilling $u_t(u_s\odot\id_t)=u_{s+t}(\id_E\odot u_{s,t})$ and $u_0$ being the canonical identification. (Also here, usually, we will leave out the mapping $u_t$ and simply identify $E\odot E_t\subset E$. The associativity condition reads, then, $(x\odot y_s)\odot z_t=x\odot(y_s\odot z_t)$.)

\brem
One may show that two \nbd{E}semigroups on the same $\sB^a(E)$ may be associated with the same product system $E^\odot$, if and only if they are conjugate by a partially isometric cocycle. While two ways to associate the same \nbd{E}semigroup with $E$ differ by a local cocycle or, what is the same, by an automorphism of $E^\odot$; see \cite[Section 7]{BhSk00} or \cite[Section 4.4]{BBLS04}.
\erem

Now we are ready to formulate the whole theorem.

\bthm\label{EE'thm}
Let $\cB\subset\sB(G)$ be a von Neumann algebra (acting nondegenerately on the Hilbert space $G$) and denote by $\cB'$ its commutant.
\begin{enumerate}
\item\label{EE'1}
The commutant establishes a one-to-one correspondence between product systems $E^\odot$ of concrete von Neumann correspondences over $\cB$ (in the sense of \cite{Ske06b}) and product systems $E'^\odot$ of concrete von Neumann correspondences over $\cB'$. Of course, $E''^\odot=E^\odot$. The product systems $E^\odot$ and $E'^\odot$ have the same morphisms $a^\odot=\bfam{a_t}_{t\in\bS}$, $a_t\in\sB^{a,bil}(E_t)=\rho'_t(\cB')'\cap\rho_t(\cB)'=\sB^{a,bil}(E'_t)$. In fact, the commutant is a bijective functor between $\cvN_\cB^\odot$ and $\cvN_{\cB'}^\odot$.

\item\label{EE'3}
Contractive units $\xi^\odot$ for $E^\odot$ correspond one-to-one to normal covariant representations $\sigma'^\odot$ of $E'^\odot$ on $G$ with $\sigma'_0=\id_{\cB'}$.

\item\label{EE'2}
Let $(E,H)\in\cvN_\cB$ and $E^\odot\in\cvN_\cB^\odot$. The ways to associate with $E^\odot$ a normal \nbd{E}semigroup $\vt$ on $\sB^a(E)$ correspond one-to-one to the normal  isometric covariant representations $\tau'^\odot$ of $E'^\odot$ on $H$ with $\tau'_0=\rho'$. Moreover,  $E$ is strongly full, if and only if $\tau'^\odot$ is faithful, and $\vt$ is an \nbd{E_0}semigroup, if and only if $\tau'^\odot$ is nondegenerate. If both is true, then necessarily every $E_t$ is strongly full, respectively, the left action of $\cB'$ on every $E'_t$ is faithful.

\item\label{EE'4}
Let $T$ be the CP-semigroup determined by either of the ingredients of a pair $(\xi^\odot,\sigma'^\odot)$ as in Number \ref{EE'3}. Then, in the sense of Number \ref{EE'2}, weak dilations $(E,\vt,\xi)$ of $T$ associated with $E^\odot$ correspond one-to-one to isometric dilations $\tau'^\odot$ of $\sigma'^\odot$. In particular, existence of one implies existence of the other.

\end{enumerate}
\ethm

\proof
\ref{EE'1}.~
This was indicated in the case of product systems of von Neumann correspondences in \cite{Ske03c}. Here we have concrete von Neumann correspondences $E_t=C_{\cB'}(\sB(G,H_t))$ and $E'_t=C_\cB(\sB(G,H_t))$ as in \cite{Ske06b}, so that a product system $E^\odot$ gives rise to a family $\bfam{E'_t}_{t\in\bS}$ of von Neumann \nbd{\cB'}correspondences. What is still missing is the product system structure of this family. Computations of this type have been detailed also in Muhly and Solel \cite{MuSo05p} (in the language of \nbd{\sigma}duals) so that here we may content ourselves with a sketchy description. The identification $E_s\odot E_t\rightarrow E_{s+t}$ is given by an operator $u_{s,t}\in\sB(E_s\odot H_t,H_{s+t})$ that intertwines both $\rho'_{s+t}$ and the canonical action $\id_{E_s}\odot\rho'_t$ of $\cB'$ on a $E_s\odot H_t$ as well as $\rho_{s+t}$ and the canonical action of $\cB$ on $E_s\odot H_t$. On the other hand, $E_s\odot H_t$ is canonically isomorphic to $E'_t\odot H_s$. Indeed, consider an element $y'_tg$ in the total subset $E'_tG$ of $H_t$. Then $x_s\odot y'_tg\mapsto y'_t\odot x_sg$ defines a unitary $E_s\odot H_t\rightarrow E'_t\odot H_s$ which intertwines the respective actions of $\cB$ and also the respective actions of $\cB'$. The image of $u_{s,t}$ as an operator $E'_t\odot H_s\rightarrow H_{s+t}$ determines an identification $E'_t\odot E'_s\rightarrow E'_{s+t}$. We leave it as an exercise to check associativity and also the statement about the morphisms.

\ref{EE'3}.~
For a single pair of correspondences $E_t\leftrightarrow E'_t$ this is more or less \cite[Lemma 2.16]{MuSo02} and the remark that follows it. (In our notations, \cite[Lemma 2.16]{MuSo02} asserts that covariant representations $(\sigma'_t,\id_{\cB'})$ of $E'_t$ on $G$ correspond one-to-one to contractions $\xi_t^*\in E_t^*\subset\sB(H_t,G)$. By the remark following \cite[Lemma 2.16]{MuSo02} normality of $\id_{\cB'}$ alone implies normality of $(\sigma'_t,\id_{\cB'})$. However, the correspondence $E_t=E_t^{**}$ has not been mentioned in \cite{MuSo02}. In the form of \nbd{\sigma}duals it appears as \cite[Theorem 3.4]{MuSo04}.) We leave it as an exercise to check that the $\xi_t$ form a unit for $E^\odot$, if and only if the $(\sigma'_t,\id_{\cB'})$ form a covariant representation of $E'^\odot$. (In the special case when $E^\odot$ is generated by $\xi^\odot$, that is if $E^\odot$ is the GNS-system of the CP-semigroup associated with $\xi^\odot$, an application of \cite[Proposition 3.1]{Ske02} shows that the statement is equivalent to \cite[Theorem 3.9]{MuSo02}. However, the product system $E^\odot=E''^\odot$, which coincides with the product system constructed in \cite{BhSk00}, has not been mentioned in \cite{MuSo02}.)

\ref{EE'2}.~
For \nbd{E_0}semigroups and nondegenerate covariant representations this is \cite[Theorem 7.4]{Ske04p} in the discrete case and \cite[Theorem 8.2]{Ske04p} in general with \cite[Remark 8.3]{Ske04p} pointing out the extensions we need here. That part of the backwards direction that constructs from a representation of $E'^\odot$ an \nbd{E} or \nbd{E_0}semigroup can be found already as a part of \cite[Theorem 3.10]{MuSo02} and goes back to \cite[Lemma 2.3]{MuSo99} for the discrete case, while the proof that the associated product system is $E^\odot$ is due to \cite[end of Section 2]{Ske03c} in the case of \nbd{E_0}semigroups and generalizes easily to \nbd{E}semigroups. We explain this very briefly. The von Neumann \nbd{\cB}module $E$ on the semigroup side and $\tau'_0$ on the covariant representation side are related by considering $\tau'_0$ as the commutant lifting associated with $E$. (From this the statement about the relation between strongly full and faithful follows.) Next if $\vt_t$ is a normal endomorphism of $\sB^a(E)$, then $\vt_t(\U)E$ factors into $E\odot E_t$ or, equivalently, $\vt_t(\U)H$ factors into $E\odot H_t\subset H$. Therefore, $\tau'_t(y'_t)\colon xg\mapsto x\odot y'_tg\in H$ defines an isometric covariant representation $(\tau'_t,\tau'_0)$ of $E'_t$ on $H$. (This representation is nondegenerate, if and only if $\vt_t$ is unital.) Similarly, if $(\tau'_t,\tau'_0)$ is an isometric covariant representation of $E'_t$ on $H$, then $\cls\tau'_t(E'_t)H$ is isomorphic to $E\odot H_t$ via $\tau'_t(y'_t)xg\mapsto x\odot y'_tg$. By $\vt_t(a)=a\odot\id_{H_t}$ we define a homomorphism $\sB^a(E)\rightarrow\sB(\cls\tau'_t(E'_t)H)\subset\sB(H)$. Actually, the range is contained in $\tau'(\cB')'=\sB^a(E)$ so that $\vt_t$ is an endomorphism. It remains to check that the semigroup property of the family $\vt$ corresponds exactly to the factorization property the family $\tau'$ must fulfill taking into account how, according to the first part, the product system structures of $E^\odot$ and of $E'^\odot$ are related. (One sees that, nicely enough, the order of elementary tensors $z'_t\odot y'_s$ acting on $xg\in H$ becomes reversed. This explains, roughly, why everything is compatible with the first part. We leave details as an exercise.)

\ref{EE'4}.~
This is simply \ref{EE'2}.~and \ref{EE'3}.~put together. Of course, a unit vector $\xi\in E\subset\sB(G,H)$ is an isometry and allows to identify $G$ as a subspace of $H$ and provides us with a projection $\xi\xi^*$ onto that subspace that compresses the isometric representation $\tau'^\odot$ to $\sigma'^\odot$. Conversely, if $G\subset H$ and the projection $p$ onto $G$ compresses $\tau'^\odot$ to $\sigma'^\odot$, then $p$ is in the intertwiner space $E$. The unit vector $\xi:=p\upharpoonright G$ has all the desired properties, that is, $\xi$ turns the \nbd{E}semigroup $\vt$ (dual to $\tau'^\odot$) into a dilation of $T$.\qed

\brem
The list of dualities may be extended. For instance, the fact that (by using \it{quasi orthonormal bases} of von Neumann \nbd{\cB}modules, as suitable substitute for orthonormal bases of Hilbert spaces) every von Neumann \nbd{\cB}module is a complemented submodule of a \it{free} von Neumann \nbd{\cB}module, may be used to prove the \it{amplification-induction theorem} on the representations $\rho'$ of $\cB'$. In the presence of invariant vector states there is a duality between CP-maps from $\cA$ to $\cB$ and CP-maps from $\cB'$ to $\cA'$ that includes a duality between \it{tensor dilations} of a CP-maps on $\cB$ and \it{extensions} from $\cB'$ to $\sB(G)$ of the dual of that CP-map; see Gohm and Skeide \cite{GoSk05}. Applying this duality to the canonical embedding of a subalgebra $\cA\subset\cB$ into $\cB$ (both in standard form) and translating back the dual map $\cB'\rightarrow\cA'$ into a map $\cB\rightarrow\cA$ via twofold \it{Tomita conjugation}, one obtains the \it{Accardi-Cecchini conditional expectation} \cite{AcCe82} that coincides with the usual conditional expectation whenever the latter exists; see also Accardi and Longo \cite{Lon84,AcLo93}.
\erem

\brem\label{coswitchrem}
The commutant switches orders in tensor products: $(E\odot F)'\cong F'\odot E'$. For von Neumann \nbd{\cB}correspondences this is \cite[Theorem 2.3]{Ske03c}. For correspondences over different von Neumann algebras (in the language of \nbd{\sigma}duals) this is \cite[Lemma 3.3]{MuSo05}. We note that it is possible to give a precise meaning to $E\odot F'$ as Connes correspondence (closely related to Sauvageot \cite{Sau80,Sau83}) and that $E\odot F'$ is canonically isomorphic to $F'\odot E$ via simply flipping elementary tensors. The isomorphisms as we fixed them in the proof of Theorem \ref{EE'thm}\eqref{EE'1} for concrete von Neumann correspondences become, then, $(E\odot F)'\cong E\odot F\cong E\odot F'\cong F'\odot E\cong F'\odot E'$ everything in the sense of Connes correspondences. Roughly speaking, in the sense of Connes correspondences the right factor in a tensor product may be replaced by its commutant. With this rule of thumb (valid for correspondences over different von Neumann algebras) we obtain a powerful tool to calculate isomorphisms in multifold tensor products of correspondences and their commutants. We will explain this and apply it systematically in \cite{Ske07p}.

We emphasize that Part \eqref{EE'2} of the theorem is also an instance of this change of order under commutant. Recall that, in the sense of Remark \ref{gencommrem} the Hilbert space $H$ with the commutant lifting $\rho'$ of $(E,H)$ is just the commutant $E'$ of $E$. In other words, we have $(E\odot E_t)'\cong E_t'\odot H$. As $u_t\colon E\odot E_t\rightarrow E$ is bilinear in the sense that $u_t$ intertwines the canonical actions of $\sB^a(E)$ on these \nbd{\cB}modules, it induces a bilinear mapping $u'_t$ between the \nbd{\rho'(\cB')}\nbd{\C}modules $E'_t\odot H$ and $H$. The representation constructed is nothing but $\tau'_t(x'_t)h=u'_t(x'_t\odot h)$. The \nbd{E}semigroup on $\sB^a(E)=\sB^{bil}(H)$ is just $\vt_t(a)=u'_t(\id'_t\odot a){u'_t}^*$.

Families like $u_t$ and $u'_t$ (in the unital, respectively, nondegenerate case) have been called left and right dilations of $E^\odot$ and of $E'^\odot$, respectively, in Skeide \cite{Ske06,Ske06p5,Ske06p6}. In \cite{Ske06} their consequent application led to a simple proof of Arveson's result \cite{Arv89,Arv90a,Arv89a,Arv90} that every Arveson system (product system of Hilbert spaces) is associated with an \nbd{E_0}semigroup. This proof has been simplified further by Arveson \cite{Arv06} and Skeide \cite{Ske06a}. A von Neumann module version with full existence results (going beyond the special Hilbert module versions \cite{Ske06p5,Ske06p6}) shall appear elsewhere \cite{Ske07p}. There we will also address questions related to continuity (or measurability) of product systems.
\erem

\section{Proof of Theorem \protect{\ref{CRdilthm}}}\label{proofSEC}

Let us discuss first the von Neumann case. We appeal to Theorem \ref{EE'thm} and the existence result Theorem \ref{CPdilthm}, where $M$ plays the role of $\cB'$ and $F^\odot$ plays the role of $E'^\odot$. For that goal we have to show that the general case boils down to the case when $\sigma_0$ is faithful and nondegenerate. But this is easy.

First of all, as $\sigma_0(\U)\sigma(y_t)\sigma_0(\U)=\sigma_t(y_t)$ we may simply pass to the subspace $\sigma_0(\U)G$ of $G$, so that now $\sigma_0$ is nondegenerate.

Then, if $\sigma_0$ is not faithful, we simply ``add'' a faithful nondegenerate representation $\wh{\sigma}_0$ of $\ker\sigma_0$ on a Hilbert space $\wh{G}$. More precisely, we pass to the covariant representation $\breve{\sigma}^\odot$ on $\breve{G}:=\wh{G}\oplus G$ that is defined by setting
\baln{
\breve{\sigma}_0
&
~:=~
\wh{\sigma}_0\oplus\sigma_0
&
\breve{\sigma}_t
&
~:=~
0\oplus\sigma_t
~~~(t>0).
}\ealn
Then every isometric dilation of $\breve{\sigma}^\odot$ compresses further to $G$, giving back $\sigma^\odot$.

Now we set $\cB:=\breve{\sigma}_0(M)'\subset\sB(\breve{G})$, so that $M\cong\cB'$ as \nbd{W^*}algebras. Put $E_t=C_{\cB'}(\sB(\breve{G},H_t))$ where $H_t:=F_t\odot\breve{G}$ with $\rho'_t$ and $\rho_t$ the canonical action of $\cB'$ and of $\cB$, repectively, so that $F^\odot\cong E'^\odot$ as product system of \nbd{W^*}correspondences. ($E^\odot$ is precisely what \cite{MuSo04,MuSo05p} would call the \it{\nbd{\breve{\sigma}_0}dual} of $F^\odot$.) Under these isomorphisms the covariant representation $\breve{\sigma}^\odot$ of $F^\odot$ induces a normal covariant representation $\sigma'^\odot$ of $E'^\odot$ on $\breve{G}$ with $\sigma'_0=\id_\cB$. So dilating as in Theorem \ref{CPdilthm} the CP-semigroup $T$ on $\cB$ associated to the unit $\xi^\odot$ corresponding to $\sigma'^\odot$ by Theorem \ref{EE'thm}\eqref{EE'3}, by Theorem \ref{EE'thm}\eqref{EE'4} we obtain an isometric dilation of $\sigma'^\odot$ and, therefore, an isometric dilation of $\sigma^\odot$.

\lf
Now let us discuss the \nbd{C^*}case, that is, $M$ is a \nbd{C^*}algebra and $F^\odot$ a product system of \nbd{C^*}correspondences over $M$ with a covariant representation $\sigma^\odot$ on a Hilbert space $G$. Our scope is simply to pass to the double commutant $F''^\odot$ of $F^\odot$ which is a product system of concrete von Neumann modules, to show that $\sigma^\odot$ extends to a normal covariant representation $\sigma''^\odot$ of that double commutant and, then, to apply the preceding discussion to $F''^\odot$. To that goal we must choose a faithful representation of $M$, and we must choose it carefully if we want that the left action of $M$ is sufficiently normal.

So let $K$ be the representation space of the \it{universal enveloping von Neumann algebra} $M^{**}$ of $M$. In this way we identify $M\subset\sB(K)$ and $M^{**}=M''$. Put $K_t:=F_t\odot K$. On $K_t$ we have the normal commutant lifting $\pi'_t$ of $M'$. (The fact that $\pi'_t$ is normal follows easily from the fact that $\sB(K,K_t)$ has enough intertwiners $y_t\in F_t$ for the action of $M'$ via $\pi'_t$ on $K_t$ and the defining action of $M'$ on $K$.) The left action of $M$ on $F_t$ gives rise to a representation $\pi_t$ of $M$ on $K_t$. By the universal property of $M''$ this representation extends to a unique normal representation $\pi''_t$ of $M''$ on $K_t$. Put $F''_t:=C_{M'}(\sB(K,K_t))\supset F_t$. Clearly, $(F''_t,K_t)$ is a concrete von Neumann correspondence over $M''$ and the product system structure of $F^\odot$ extends uniquely to $F''^\odot$.

We return to the covariant representation $\sigma^\odot$ of $F^\odot$. First, we observe that $\sigma_0$ extends uniquely to a normal representation $\sigma''_0$ of $M''$. Then, as in the proof Theorem \ref{EE'thm}\eqref{EE'3}, $\sigma^\odot$ gives rise to a unit $\zeta'^\odot$ for $F'^\odot:=(F'')'^\odot$. Further, each $\zeta'_t\in F'_t$ gives rise to a normal representation $(\sigma''_t,\sigma''_0)$ of $F''_t$, which clearly extends $(\sigma_t,\sigma_0)$. (This is exactly the statement of the remark following \cite[Theorem 2.16]{MuSo02}.) Of course, the $\sigma''_t$ form a normal representation $\sigma''^\odot$ of $F''^\odot$ so that now we are ready to apply the \nbd{W^*}version of Theorem \ref{CRdilthm}, obtaining an isometric dilation $\tau''^\odot$ of $\sigma''^\odot$ that restricts to an isometric dilation $\tau^\odot$ of $\sigma^\odot$.

\brem
We mention that if the weak dilation that we use as input for Theorem \ref{EE'thm}\eqref{EE'4} is the unique one from Theorem \ref{CPdilthm} fulfilling the Condition \ref{CPd2}, then the corresponding isometric dilation is also \hl{minimal} in the sense of \cite[Theorem 3.7]{MuSo02} (that is, the smallest subspace of $H$ invariant for the isometric dilation is $H$) and determined uniquely by this condition.
\erem

\vfill

\setlength{\baselineskip}{2.5ex}


\newcommand{\Swap}[2]{#2#1}\newcommand{\Sort}[1]{}
\providecommand{\bysame}{\leavevmode\hbox to3em{\hrulefill}\thinspace}
\providecommand{\MR}{\relax\ifhmode\unskip\space\fi MR }
\providecommand{\MRhref}[2]{%
  \href{http://www.ams.org/mathscinet-getitem?mr=#1}{#2}
}
\providecommand{\href}[2]{#2}


\end{document}